\documentclass{llncs}
\usepackage[utf8]{inputenc}
\usepackage{tikz}
\usetikzlibrary{tikzmark}

\usepackage{float}
\usepackage{caption}
\usepackage{subcaption}

\usepackage[numbers]{natbib}
\begin{document}

\title{An application of neural networks to a problem in knot theory and group theory (untangling braids)}

\titlerunning{An application of neural networks to untangling braids}

\author{Alexei Lisitsa\inst{1}
\and
Mateo Salles\inst{2} \and
Alexei Vernitski\inst{2}
}
\authorrunning{A. Lisitsa, M. Salles, A. Vernitski}
%
\institute{University of Liverpool,
  Liverpool, UK \and
University of Essex, UK
}

\maketitle

\begin{abstract}
We report on our success on solving the problem of untangling braids up to length 20 and width 4. We use feed-forward neural networks in the framework of reinforcement learning to train the agent to choose Reidemeister moves to untangle braids in the minimal number of moves.

\end{abstract}

\subsubsection{Context and novelty}

Recent studies \cite{gukov2021learning,kauffman2020rectangular} report that they use elements of deep learning to untangle braids (or knots), but do not provide details. Here we present all details of our implementation. Also, we have published our code, and it can be accessed on GitHub \cite{ourcode}.

In our previous research \cite{khan2021untangling} we used Q-learning to untangle braids, but in this study, as we train the agent to untangle larger braids (up to 20 crossings on up to 4 strands), we combine reinforcement learning and deep learning.

When we generate random trivial braids for training and testing, we use a version of a classical group-theoretical algorithm for checking if a braid is trivial \cite[Section 1.5]{kassel2008braid}. For the number of crossings 16--20 we observed that our neural network untangles braids faster than the group-theoretical algorithm proves that they are trivial.


\subsubsection{Braids}

\begin{figure}[h]
\begin{tikzpicture}[scale=0.50]
\draw (0,0) -- (2,0);
\draw (0,4) .. controls (1,4) and (1,2) .. (2,2);
\draw (0,2) .. controls (0.4,2)  .. (0.9,2.85);
\draw (2,4) .. controls (1.6,4)  .. (1.1,3.15);
\draw (0,6) -- (2,6);
\draw (2,0) -- (4,0);
\draw (2,2) -- (4,2);
\draw (2,6) .. controls (3,6) and (3,4) .. (4,4);
\draw (2,4) .. controls (2.4,4)  .. (2.9,4.85);
\draw (4,6) .. controls (3.6,6)  .. (3.1,5.15);
\draw (4,0) -- (6,0);
\draw (4,2) -- (6,2);
\draw (4,6) .. controls (5,6) and (5,4) .. (6,4);
\draw (4,4) .. controls (4.4,4)  .. (4.9,4.85);
\draw (6,6) .. controls (5.6,6)  .. (5.1,5.15);
\draw (6,0) -- (8,0);
\draw (6,2) -- (8,2);
\draw (6,4) .. controls (7,4) and (7,6) .. (8,6);
\draw (6,6) .. controls (6.4,6)  .. (6.9,5.15);
\draw (8,4) .. controls (7.6,4)  .. (7.1,4.85);
\draw (8,0) -- (10,0);
\draw (8,2) -- (10,2);
\draw (8,4) .. controls (9,4) and (9,6) .. (10,6);
\draw (8,6) .. controls (8.4,6)  .. (8.9,5.15);
\draw (10,4) .. controls (9.6,4)  .. (9.1,4.85);
\draw (10,2) .. controls (11,2) and (11,0) .. (12,0);
\draw (10,0) .. controls (10.4,0)  .. (10.9,0.85);
\draw (12,2) .. controls (11.6,2)  .. (11.1,1.15);
\draw (10,4) -- (12,4);
\draw (10,6) -- (12,6);
\draw (12,0) -- (14,0);
\draw (12,2) -- (14,2);
\draw (12,6) .. controls (13,6) and (13,4) .. (14,4);
\draw (12,4) .. controls (12.4,4)  .. (12.9,4.85);
\draw (14,6) .. controls (13.6,6)  .. (13.1,5.15);
\draw (14,0) .. controls (15,0) and (15,2) .. (16,2);
\draw (14,2) .. controls (14.4,2)  .. (14.9,1.15);
\draw (16,0) .. controls (15.6,0)  .. (15.1,0.85);
\draw (14,4) -- (16,4);
\draw (14,6) -- (16,6);
\draw (16,0) -- (18,0);
\draw (16,2) .. controls (17,2) and (17,4) .. (18,4);
\draw (16,4) .. controls (16.4,4)  .. (16.9,3.15);
\draw (18,2) .. controls (17.6,2)  .. (17.1,2.85);
\draw (16,6) -- (18,6);
\draw (18,0) -- (20,0);
\draw (18,2) -- (20,2);
\draw (18,4) .. controls (19,4) and (19,6) .. (20,6);
\draw (18,6) .. controls (18.4,6)  .. (18.9,5.15);
\draw (20,4) .. controls (19.6,4)  .. (19.1,4.85);
\draw (20,0) -- (22,0);
\draw (20,2) .. controls (21,2) and (21,4) .. (22,4);
\draw (20,4) .. controls (20.4,4)  .. (20.9,3.15);
\draw (22,2) .. controls (21.6,2)  .. (21.1,2.85);
\draw (20,6) -- (22,6);
\draw (22,0) -- (24,0);
\draw (22,2) -- (24,2);
\draw (22,6) .. controls (23,6) and (23,4) .. (24,4);
\draw (22,4) .. controls (22.4,4)  .. (22.9,4.85);
\draw (24,6) .. controls (23.6,6)  .. (23.1,5.15);\end{tikzpicture}
\caption{A braid on $4$ strands}
\label{fig:braid}
\end{figure}
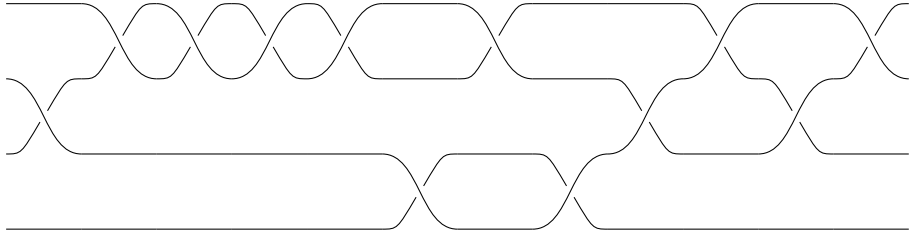

Braids are mathematical objects from low-dimensional topology or, to be more precise, knot theory (that is, the study of the relative position of curves in the space). A \emph{braid} on $n$ strands consists of $n$ ropes whose left-hand ends are fixed one under another and whose right-hand ends are fixed one under another; you can imagine that the braid is laid out on a table, and the ends of the ropes are attached to the table with nails. Figure \ref{fig:braid} shows an example of braids on $4$ strands with $12$ crossings; this is an example of a braid which our code successfully untangles. Braids are important because, on the one hand, they are useful building blocks of knots and other constructions of low-dimensional topology and, on the other hand, have a simple structure and can be conveniently studied using mathematics and, as in this study, experimented with using computers.

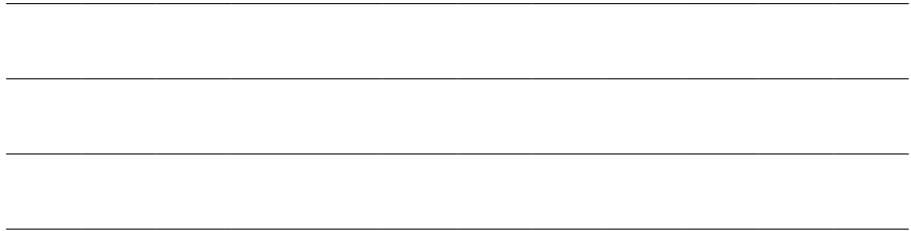
\begin{figure}[h]
\begin{tikzpicture}[scale=0.50]
\draw (0,0) -- (2,0);
\draw (0,2) -- (2,2);
\draw (0,4) -- (2,4);
\draw (0,6) -- (2,6);
\draw (2,0) -- (4,0);
\draw (2,2) -- (4,2);
\draw (2,4) -- (4,4);
\draw (2,6) -- (4,6);
\draw (4,0) -- (6,0);
\draw (4,2) -- (6,2);
\draw (4,4) -- (6,4);
\draw (4,6) -- (6,6);
\draw (6,0) -- (8,0);
\draw (6,2) -- (8,2);
\draw (6,4) -- (8,4);
\draw (6,6) -- (8,6);
\draw (8,0) -- (10,0);
\draw (8,2) -- (10,2);
\draw (8,4) -- (10,4);
\draw (8,6) -- (10,6);
\draw (10,0) -- (12,0);
\draw (10,2) -- (12,2);
\draw (10,4) -- (12,4);
\draw (10,6) -- (12,6);
\draw (12,0) -- (14,0);
\draw (12,2) -- (14,2);
\draw (12,4) -- (14,4);
\draw (12,6) -- (14,6);
\draw (14,0) -- (16,0);
\draw (14,2) -- (16,2);
\draw (14,4) -- (16,4);
\draw (14,6) -- (16,6);
\draw (16,0) -- (18,0);
\draw (16,2) -- (18,2);
\draw (16,4) -- (18,4);
\draw (16,6) -- (18,6);
\draw (18,0) -- (20,0);
\draw (18,2) -- (20,2);
\draw (18,4) -- (20,4);
\draw (18,6) -- (20,6);
\draw (20,0) -- (22,0);
\draw (20,2) -- (22,2);
\draw (20,4) -- (22,4);
\draw (20,6) -- (22,6);
\draw (22,0) -- (24,0);
\draw (22,2) -- (24,2);
\draw (22,4) -- (24,4);
\draw (22,6) -- (24,6);\end{tikzpicture}
\caption{The trivial braid diagram on $4$ strands}
\label{fig:trivial}
\end{figure}

The braid in Figure \ref{fig:braid} can be untangled, that is, all crossings can be removed by moving certain parts of strands up or down, as needed (without touching the ends of the ropes); after the braid is untangled, the braid diagram will look as in Figure \ref{fig:trivial}, which shows the \emph{trivial braid diagram}. Not every braid can be untangled. Those braids that can be untangled are called \emph{trivial} braids. The task that we successfully solve in this research is training a neural network to untangle braids. 

\subsubsection{Reidemeister moves}

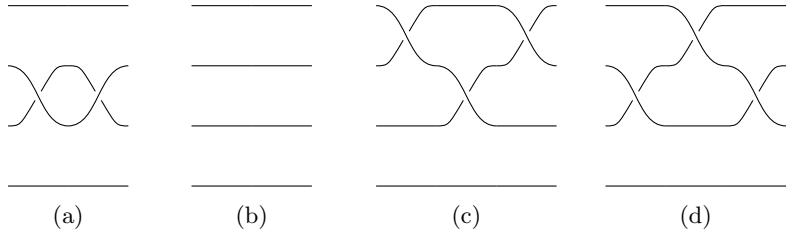
\begin{figure}
\centering
\begin{subfigure}{.2\textwidth}
  \centering
\begin{tikzpicture}[scale=0.40]
\draw (0,0) -- (2,0);
\draw (0,4) .. controls (1,4) and (1,2) .. (2,2);
\draw (0,2) .. controls (0.4,2)  .. (0.9,2.85);
\draw (2,4) .. controls (1.6,4)  .. (1.1,3.15);
\draw (0,6) -- (2,6);
\draw (2,0) -- (4,0);
\draw (2,2) .. controls (3,2) and (3,4) .. (4,4);
\draw (2,4) .. controls (2.4,4)  .. (2.9,3.15);
\draw (4,2) .. controls (3.6,2)  .. (3.1,2.85);
\draw (2,6) -- (4,6);
\end{tikzpicture}  \caption{}
  \label{fig:sub1}
\end{subfigure}%
\begin{subfigure}{.2\textwidth}
  \centering
\begin{tikzpicture}[scale=0.40]
\draw (0,0) -- (2,0);
\draw (0,2) -- (2,2);
\draw (0,4) -- (2,4);
\draw (0,6) -- (2,6);
\draw (2,0) -- (4,0);
\draw (2,2) -- (4,2);
\draw (2,4) -- (4,4);
\draw (2,6) -- (4,6);
\end{tikzpicture}  \caption{}
  \label{fig:sub2}
\end{subfigure}
\begin{subfigure}{.25\textwidth}
  \centering
\begin{tikzpicture}[scale=0.40]
\draw (0,0) -- (2,0);
\draw (0,2) -- (2,2);
\draw (0,6) .. controls (1,6) and (1,4) .. (2,4);
\draw (0,4) .. controls (0.4,4)  .. (0.9,4.85);
\draw (2,6) .. controls (1.6,6)  .. (1.1,5.15);
\draw (2,0) -- (4,0);
\draw (2,4) .. controls (3,4) and (3,2) .. (4,2);
\draw (2,2) .. controls (2.4,2)  .. (2.9,2.85);
\draw (4,4) .. controls (3.6,4)  .. (3.1,3.15);
\draw (2,6) -- (4,6);
\draw (4,0) -- (6,0);
\draw (4,2) -- (6,2);
\draw (4,6) .. controls (5,6) and (5,4) .. (6,4);
\draw (4,4) .. controls (4.4,4)  .. (4.9,4.85);
\draw (6,6) .. controls (5.6,6)  .. (5.1,5.15);
\end{tikzpicture}  \caption{}
  \label{fig:sub3}
\end{subfigure}%
\begin{subfigure}{.25\textwidth}
  \centering
\begin{tikzpicture}[scale=0.40]
\draw (0,0) -- (2,0);
\draw (0,4) .. controls (1,4) and (1,2) .. (2,2);
\draw (0,2) .. controls (0.4,2)  .. (0.9,2.85);
\draw (2,4) .. controls (1.6,4)  .. (1.1,3.15);
\draw (0,6) -- (2,6);
\draw (2,0) -- (4,0);
\draw (2,2) -- (4,2);
\draw (2,6) .. controls (3,6) and (3,4) .. (4,4);
\draw (2,4) .. controls (2.4,4)  .. (2.9,4.85);
\draw (4,6) .. controls (3.6,6)  .. (3.1,5.15);
\draw (4,0) -- (6,0);
\draw (4,4) .. controls (5,4) and (5,2) .. (6,2);
\draw (4,2) .. controls (4.4,2)  .. (4.9,2.85);
\draw (6,4) .. controls (5.6,4)  .. (5.1,3.15);
\draw (4,6) -- (6,6);
\end{tikzpicture}  \caption{}
  \label{fig:sub4}
\end{subfigure}
\caption{Examples of Reidemeister moves}
\label{fig:Reidemeister}
\end{figure}

When one studies braids (or knots) and how to untangle them, the untangling process is split into elementary local changes, affecting $2$ or $3$ consecutive crossings, called Reidemeister moves \cite{kassel2008braid}. Transition from (a) to (b) in Figure \ref{fig:Reidemeister} is one example of a Reidemeister move. In this example, the Reidemeister move removes two crossings from a braid; thus, if our aim is to untangle the braid, it seems like a good move to use. However, not every Reidemeister move removes crossings from a braid; some Reidemeister moves only prepare groundwork for removing crossings; one example is transition from (c) to (d) in Figure \ref{fig:Reidemeister}. Thus, the challenge for artificial intelligence is to learn which Reidemeister moves it should use to untangle a braid, even though initially it might be not clear how these moves contribute to untangling the braid.

\subsubsection{Reinforcement learning}

The way in which we use the neural network is to produce a recommendation for the agent as to what Reidemeister move can be applied to the current configuration of the braid. However, this move is only one move in a sequence of moves which untangles the braid given to the agent. Thus, we use a neural network in the framework of reinforcement learning. 

We give the agent a small negative reward whenever it performs a Reidemeister move, and a large positive reward when the agent has successfully produced the trivial braid diagram (following the example of agents finding a way out of a maze \cite{sutton2018reinforcement}). Thus, the agent is encouraged to untangle the braid, and to untangle the braid in the minimal number of moves.

\subsubsection{The input of the neural network}

Considering a braid from left to right, we record a clockwise (or anti-clockwise) crossing of strands in positions $i$ and $i+1$ (counting from the bottom) as the number $i$ (or $-i$). Thus, the braid in Figure \ref{fig:braid} is encoded as the list of numbers $2, 3, 3, -3, -3, 1, 3, -1, -2, -3, -2, 3$. We were prepared to transform this encoding into one-hot encoding, if needed, but we were successful with this simple encoding. This list of numbers is fed into a neural network.

For our experiments, we fix the length of the braid, as we did in our previous research \cite{khan2021untangling}. If there is no intersection of strands in a certain part of a braid, this part of the braid is denoted by $0$. Thus, in the context of braids with up to $12$ crossings, the trivial braid diagram (shown in Figure \ref{fig:trivial}) is encoded by the list of $12$ $0$s.

Not every type of Reidemeister move can be applied at every position in a given braid; for example, the transition from (a) to (b) in Figure \ref{fig:Reidemeister} cannot be applied anywhere in the braid in Figure \ref{fig:braid} because the fragment (a) does not match any part of the braid in Figure \ref{fig:braid}. A binary list indicating at which positions which Reidemeister moves can be used in the current braid is also fed into the neural network (without giving the neural network any indication as to how the values in this list can be interpreted). 

\subsubsection{The output of the neural network}

We use a softmax output layer, in which nodes correspond to all positions where Reidemeister moves can be applied in a braid. Not all Reidemeister moves can be applied to the current braid in all positions; the nodes in the output layer corresponding to impossible moves are ignored (as it is done in AI agents playing board games \cite{silver2017mastering}). Out of those nodes which correspond to possible moves, the maximal one is chosen.

\subsubsection{Network topology, training and performance}

We use a feed-forward neural network with one hidden layer, whose size is the number of inputs times the number of strands. We were prepared to use a more complicated topology, but we were successful with this simple choice.

We run 10,000 training iterations, with 2 episodes per iteration, and each episode has a maximum of 50 steps. In each episode, we draw braid out of a randomly generated trivial braid data set. Braids with less crossings lead to better accuracy as there are less moves needed.

The neural network gets trained almost immediately and shows accuracy of about $93\%$ on the braids with $10$ crossings on $4$ strands. Larger sizes of braids also train well, but the performance strongly fluctuates depending on the maximum number of steps allowed to untangle the braids. Having more steps available to untangle a braid translates to an improvement in the accuracy rate.

\begin{figure}[H]
    \centering
    \subfloat[]{{\includegraphics[scale=0.44]{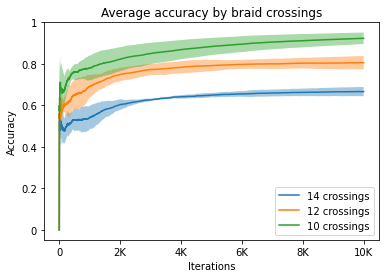} }}
    \subfloat[]{{\includegraphics[scale=0.44]{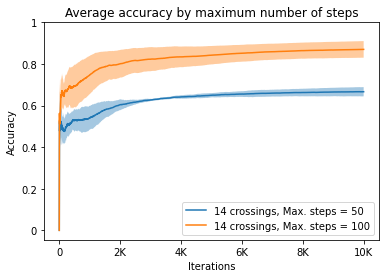} }}
    \caption{Performance for different number of crossings (a), and maximum number of steps (b)}
    \label{fig:8}
\end{figure}

\subsubsection{References}

\bibliography{main}

\end{document}